%% file: RP2.tex
\documentclass[reqno]{gokova}
\usepackage{epsfig,graphicx}

\begin{document}
\input goksty.tex      

\setcounter{page}{1}
\volume{13}

\title[]{Twisting $4$-manifolds along ${\R\P}^2$}
\author[]{Selman Akbulut}

\thanks{The author is partially supported by NSF grant DMS 0905917}

\address{Department  of Mathematics, Michigan State University,  MI, 48824}
\email{akbulut@math.msu.edu}


\begin{abstract}
We prove that the Dolgachev surface $E(1)_{2,3}$ (which is an exotic copy of the elliptic surface
 $E(1)={\C\P}^{2} \# 9 \bar{\C\P}^2$) can be obtained from $E(1)$ by twisting along a simple ``plug'', in particular it can be obtained  from $E(1)$ by twisting along ${\R\P}^2$.
\end{abstract}

\maketitle

\section{Introduction}

Given a smooth $4$-manifold $M^4$, what is the minimal genus $g$ of an imbedded surface $\Sigma_{g}\subset M^4$, such that twisting $M$ along $\Sigma $  produces an exotic copy of $M$?
Here twisting means cutting out a tubular neighborhood of $\Sigma $ and regluing back by a nontrivial diffeomorphism.  When $g>1$ we don't get anything new (bacause by \cite{o} pp.133 \footnote{we thank Cameron Gordon for pointing out this reference} any diffeomorphism of a circle bundle over $\Sigma_{g}$ can be isotoped to preserve the fiber, and hence it extends to the corresponding disk bundle).
The case $g=1$ is the well  known``logarithmic transform" operation, which can change the smooth structure in some cases; in fact the first example of a closed exotic manifold  found by Donaldson \cite{d} was the Dolgachev surface $E(1)_{2,3}$ which is obtained from $E(1) ={\C\P}^2\# 9\bar{\C\P}^2$ by two log transforms . The $g=0$ case is not well understood, twisting along $S^2$ is usually  called ``Gluck construction" and we don't know if this operation changes the smooth structure of an any orientable manifold, but  there is an example of non-orientable manifold  which the Gluck construction changes its smooth structure \cite{a1}. The interesting case of $\Sigma  = {\R\P}^2$ was studied  indirectly in \cite{ay1} under the guise of {\it plugs}, which are more general objects.  Recall that Figure 1 describes the tubular neighborhood $W$ of ${\R\P}^2$ in $S^4$ as a disc bundle over ${\R\P}^2$ (e.g. \cite{a2}):

 \begin{figure}[ht]  \begin{center}  
\includegraphics[width=.2\textwidth]{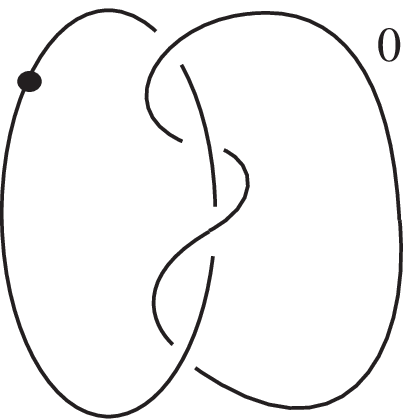}   
\caption{W} 
\end{center}
\end{figure} 

If we attach a $2$-handle to $W$ as in Figure 2 we obtain an interesting manifold, which is the $W_{1,2} $  ``plug''  of \cite{ay1}. Recall  \cite{ay1}, a {\it  plug}  $(P, f)$  of $M^4$ is a codimension zero Stein submanifold $P\subset M$ with an involution $f : \partial P \to \partial P$, such that $f$ does not extend to a homemorphism inside; and the operation $N\cup_{id}P\mapsto N\cup_{f}P$ of removing $P$ from $M$ and regluing it to its complement $N$ by $f$,  changes the smooth structure of $M$
(this operation is called  a ``{\it plug twisting}''). For example  the involution $f: \partial W_{1,2}\to \partial W_{1,2}$ is induced from $180^{0}$ rotation of the Figure 2 , e.g. it maps the (red and blue)  loops to each other $\alpha \leftrightarrow \beta$.

\begin{figure}[ht]  \begin{center}  
\includegraphics[width=.3\textwidth]{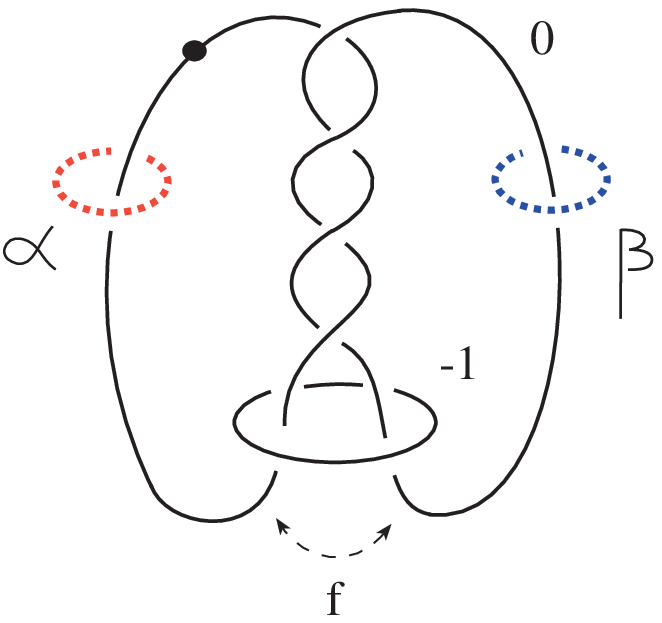}   
\caption{$W_{1,2}$} 
\end{center}
\end{figure} 

\noindent Notice that  the twisting along $W_{1,2}$ is induced by twisting along ${\R\P}^2$ inside (i.e. cutting out $W$ and regluing by the involution induced by the rotation). In  \cite{ay1} some examples of changing smooth structures via plug twisting were given, including twisting the $W_{1,2}$ plug.  Here we prove that  by twisting along a $W_{1,2}$ plug (in particular twisting along ${\R\P}^2$)  we can completely decompose  the Dolgachev surface $E(1)_{2,3}$. The following theorem should be considered as a structure theorem for the Dolgachev surface complementing Theorem 1 of \cite{a3}, where it was shown that a  ``{\it cork twisting}''  also completely decomposes $E(1)_{2,3}$.

\begin{thm} $E(1)_{2,3}$ is obtained by plug twisting $E(1)$ along  $W_{1,2}$, i.e. we can decompose   $E(1)=N\cup_{id}W_{1,2}$, so that $E(1)_{2,3}= N\cup_{f}W_{1,2}$.
\end{thm}

\proof By cancelling the $1$- and $2$-handle pair of Figure 2 we obtain Figure 3, which is an alternative picture of $W_{1,2}$. By inspecting the  diffeomorphism Figure 2 $\mapsto $ Figure 3 we see that  the involution  $f$ twists the tubular neighborhood of $\alpha$ once, while mapping to $\beta $.

 \begin{figure}[ht]  \begin{center}  
\includegraphics[width=.30\textwidth]{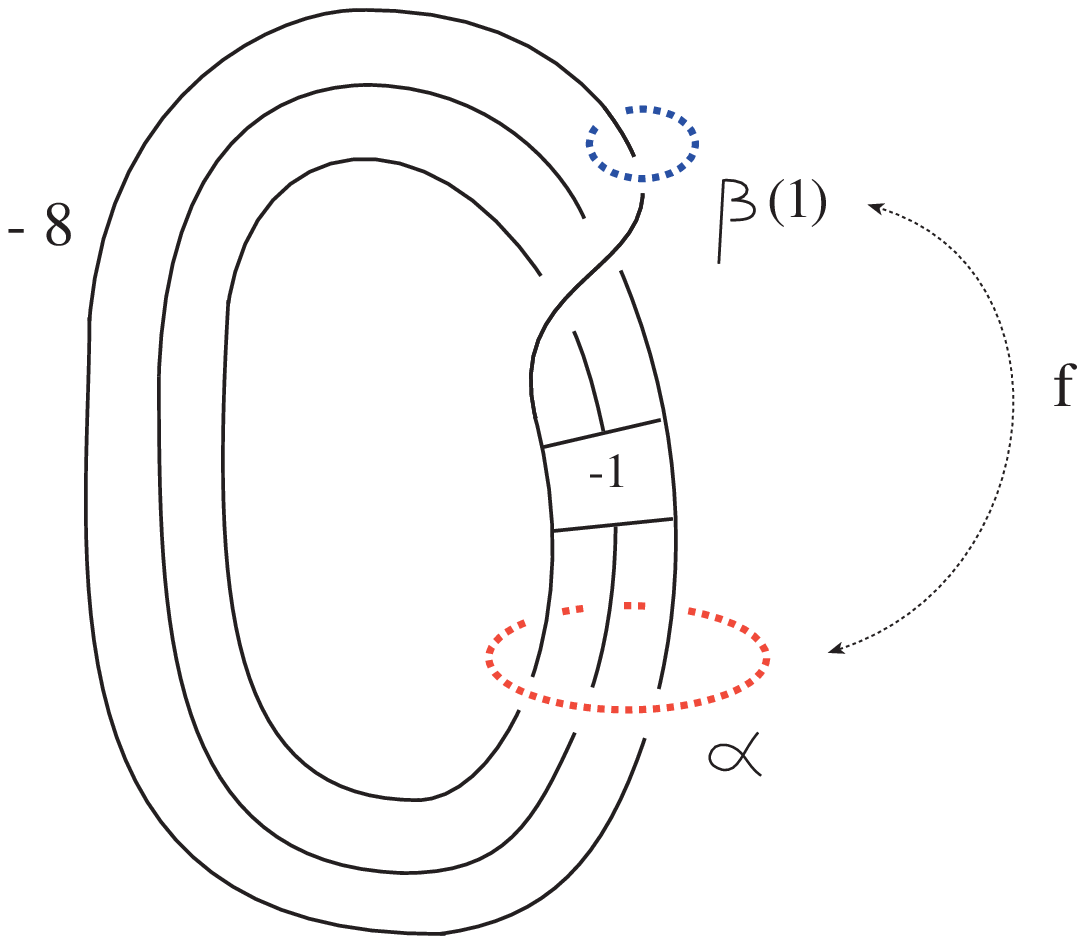}   
\caption{$W_{1,2}$} 
\end{center}
\end{figure} 

By attaching a chain of eight $2$-handles to $-W_{1,2}$ (the mirror image of Figure 3) and a $+1$ framed $2$-handle to  $\alpha$, we obtain Figure 4, which is a handlebody of $E(1)$ given in \cite{a3}. In Figure 4 performing $W_{1,2}$ plug twist to E(1) has the effect of replacing the $+1$-framed 2-handle attached to $\alpha $, with a zero framed $2$-handle attached to $\beta$. Here  the complement of $W_{1,2}$ in $E(1)$ is  the submanifold $N$  consisting of the zero framed 2-handle (the cusp) and the chain of eight $2$-handles, and the plug twisting is the operation:
$N\cup \alpha^{+1} \mapsto N\cup \beta^{0}$ (as seen from $N$).

 \begin{figure}[ht]  \begin{center}  
\includegraphics[width=.47\textwidth]{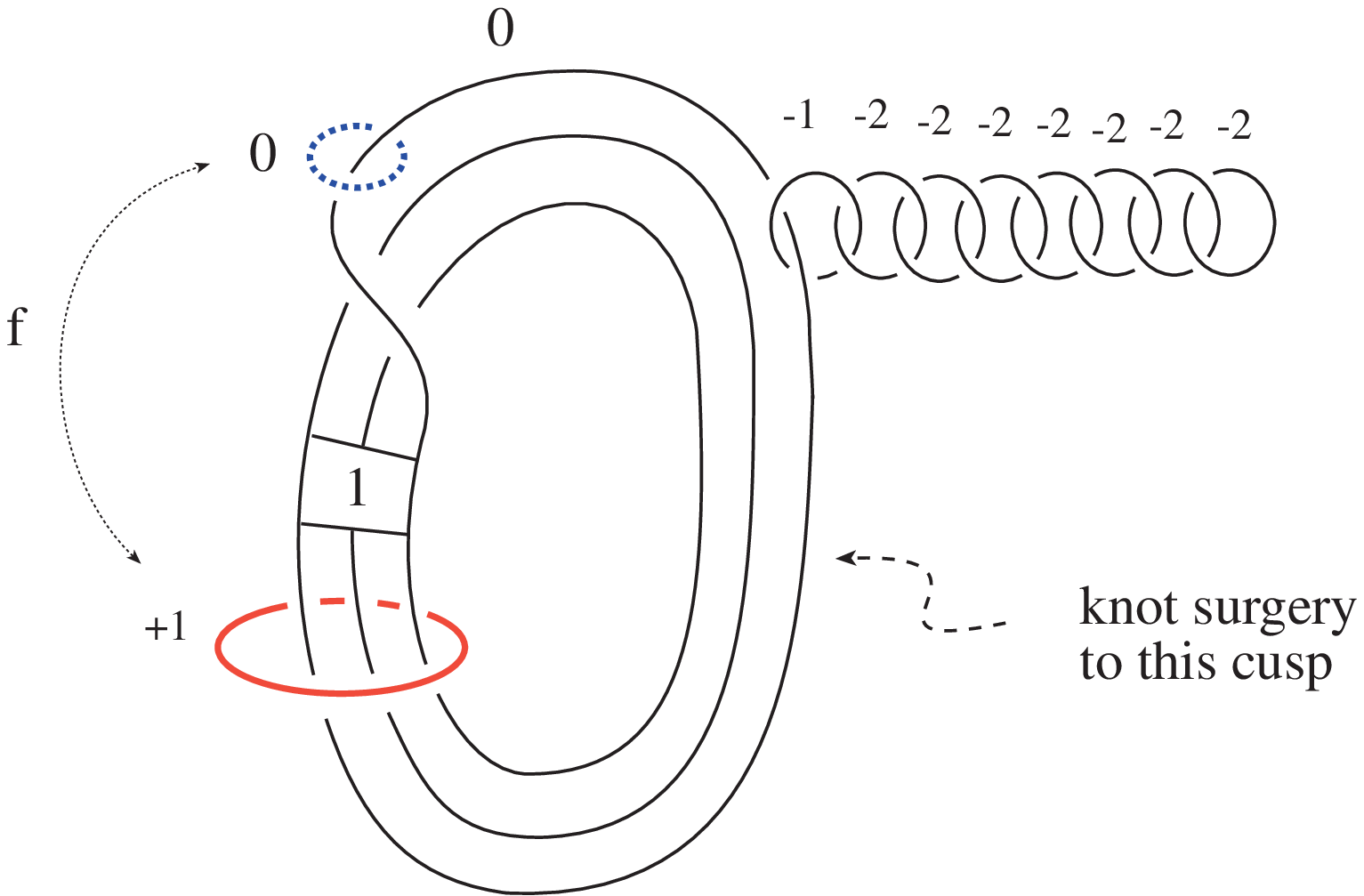}   
\caption{E(1)} 
\end{center}
\end{figure}

Therefore  the  plug twisting  of $E(1)$ along  $W_{1,2}$ gives  Figure 5. After sliding over $\beta $, the chain of eight $2$-handles become free from the rest of the figure, giving a splitting: 
$Q \#8 \bar{\C\P}^2$, where $Q$ is the cusp with the trivially linking zero framed cicle, hence $Q=S^2\times S^2$. So the Figure 5 is just $S^2\times S^2 \# 8 \bar{\C\P}^2 =E(1)$. 

 \begin{figure}[ht]  \begin{center}  
\includegraphics[width=.47\textwidth]{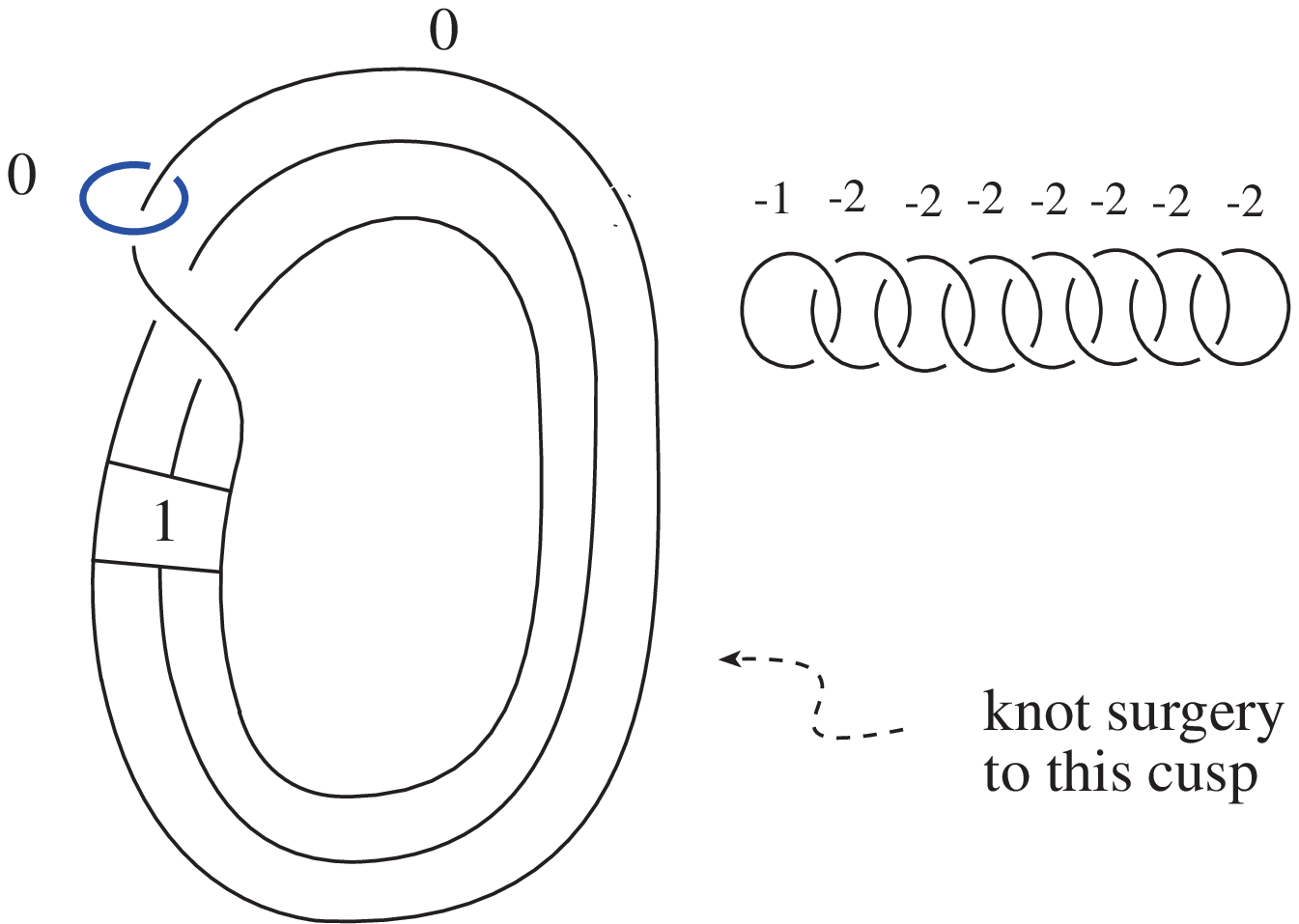}   
\caption{} 
\end{center}
\end{figure} 

Next notice that if we first perform  a ``knot surgery''  operation  $E(1) \mapsto E(1)_{K}$ by a knot $K$,  along the cusp inside of Figure 4, and then do the plug twist along $W_{1,2}$ (notice the  cusp is disjoint from the plug since it lies in $N$) we get the similar splitting except this time resulting: $Q_{K}  \# 8 \bar{\bf CP}^2$, where  $Q_{K}$ is the knot surgered $Q$. Notice  the manifold  $Q=S^2\times S^2$ is obtained by doubling the cusp, and $Q_{K}$ is obtained  by doing knot surgery to one of these cusps. In  Theorem 4.1 of  \cite{a4} it was shown that when $K$ is the trefoil knot  then $Q_{K}=S^2 \times S^2$. Also recall that when $K$ is the trefoil knot we have the identification with the Dolgachev surface $E(1)_{K}=E(1)_{2,3}$  (e.g. \cite{a3}). \qed

\vspace{.1in}

 \begin{rem} If we could identify $Q_{K}$ with $S^2\times S^2$  for infinitely many knots $K$ with distinct Alexander polynomials, we would have infinitely many transforms $E(1)\mapsto E(1)_{K}$ obtained by plug twistings along $W_{1,2}$. This would give infinitely many non-isotopic imbeddings   $W_{1,2}\subset E(1)$, similar to the examples  in \cite{ay2}. In the absence of such identification we can only conclude  that $W_{1,2}$ is a plug of infinitely many distinct exotic copies $E(1)_{K} $ of $E(1)$. 
 \end{rem}

\begin{rem} Recall that $\partial W$ is the quaternionic 3-manifold, which is the quotient  of  $S^3$ by the free action of the quaternionic group of order eight $G =< i, j, k \;|\;  i^2 = j^2 = k^2 = -1, ij = k, jk = i, ki = j >$ (e.g. \cite{a2}). This manifold  is a positively curved space-form and an L space (Floer homology groups vanish).  Hence  the change of smooth structure of $E(1)$ by twisting $W$ is due to the change of $Spin^{c}$ structures, rather than permuting the Floer homology by the involution as in \cite{a3}, \cite{ad}. 
\end{rem}

\end{document}

%% file: goksty.tex
\def\E{\ifmmode{\mathbb E}\else{$\mathbb E$}\fi} 
\def\N{\ifmmode{\mathbb N}\else{$\mathbb N$}\fi} 
\def\R{\ifmmode{\mathbb R}\else{$\mathbb R$}\fi} 
\def\Q{\ifmmode{\mathbb Q}\else{$\mathbb Q$}\fi} 
\def\C{\ifmmode{\mathbb C}\else{$\mathbb C$}\fi} 
\def\H{\ifmmode{\mathbb H}\else{$\mathbb H$}\fi} 
\def\Z{\ifmmode{\mathbb Z}\else{$\mathbb Z$}\fi} 
\def\P{\ifmmode{\mathbb P}\else{$\mathbb P$}\fi} 
\def\T{\ifmmode{\mathbb T}\else{$\mathbb T$}\fi} 
\def\SS{\ifmmode{\mathbb S}\else{$\mathbb S$}\fi} 
\def\DD{\ifmmode{\mathbb D}\else{$\mathbb D$}\fi} 

\renewcommand{\a}{\alpha}
\renewcommand{\b}{\beta}
\renewcommand{\d}{\delta}
\newcommand{\D}{\Delta}
\newcommand{\e}{\varepsilon}
\newcommand{\g}{\gamma}
\newcommand{\G}{\Gamma}
\newcommand{\la}{\lambda}
\newcommand{\La}{\Lambda}
\newcommand{\n}{\nabla}
\newcommand{\var}{\varphi}
\newcommand{\s}{\sigma}
\newcommand{\Sig}{\Sigma}
\renewcommand{\t}{\tau}
\renewcommand{\th}{\theta}
\renewcommand{\O}{\Omega}
\renewcommand{\o}{\omega}
\newcommand{\z}{\zeta}
\newcommand{\B}{{\mathbb B}}

\newcommand{\ben}{\begin{enumerate}}
\newcommand{\een}{\end{enumerate}}
\newcommand{\be}{\begin{equation}}
\newcommand{\ee}{\end{equation}}
\newcommand{\bea}{\begin{eqnarray}}
\newcommand{\eea}{\end{eqnarray}}
\newcommand{\bc}{\begin{center}}
\newcommand{\ec}{\end{center}}

\newtheorem{thm}{Theorem}[section]
\newtheorem{cor}[thm]{Corollary}
\newtheorem{lem}[thm]{Lemma}
\newtheorem{prop}[thm]{Proposition}
\newtheorem{ax}{Axiom}
\newtheorem{conj}[thm]{Conjecture}

\theoremstyle{definition}
\newtheorem{defn}{Definition}[section]

\theoremstyle{remark}
\newtheorem{rem}{\rm\bfseries{Remark}}[section]
\newtheorem*{notation}{Notation}

\newtheorem{ques}{\rm\bfseries{Question}}[section]
\newtheorem{cons}[rem]{\rm\bfseries{Construction}}
\newtheorem{exm}[rem]{\rm\bfseries{Example}}
